\documentclass[12pt]{amsart}
\usepackage{mathrsfs}
\usepackage{amsgen}
\usepackage{amscd}
\usepackage{latexsym}
\usepackage{amsthm}
\usepackage{color}
\usepackage{enumerate}
\usepackage{enumitem}
\usepackage{amssymb,amsmath,graphicx,amsfonts,euscript}
\usepackage{color}
\usepackage{empheq}
\usepackage{cases}
\usepackage{comment}
\usepackage{hyperref}
\usepackage{booktabs}
\usepackage{todonotes}
\usepackage{setspace} 
\usepackage{multirow}
\usepackage[numbers,sort]{natbib}
\usepackage{titlesec}

\makeatletter
\renewcommand{\section}{%
  \@startsection{section}{1}{\z@}%
  {-3.5ex \@plus -1ex \@minus -.2ex}%
  {2.3ex \@plus .2ex}%
  {\normalfont\Large\bfseries}}

\usepackage[margin=3cm, a4paper]{geometry}
\titlespacing*{\subsection}{0pt}{2ex}{0.5ex}

 \def\norm#1{\|#1\|}

\def\H{{\cal H}}

\def\H1{H^1(\R)}

\numberwithin{equation}{section}
\newtheorem{theorem}{Theorem}[section]

\theoremstyle{definition}

\newtheorem*{exam*}{Examples}
\newtheorem{example}{Example}[section]


\begin{document}

\setcounter{page}{1}

\title[Numerical studies of Fully fractional PDEs...]{Numerical Analysis of Fully Nonlocal Partial Differential Equations Using Graded Meshes}

\author[]{Tapaswini Sahu}
\address{Department of Mathematical Sciences}
\email{tapaswinis@iitbhilai.ac.in}
\thanks{}

\author[]{Swati Yadav}
\address{Department of Mathematical Sciences}
\email{yswati.iitbhu@gmail.com}
\thanks{}

\subjclass[2020]{35R11; 65M06; 65M12}
\keywords{Fractional diffusion equation, Caputo derivative, $L1$ scheme, graded mesh, fractional Laplacian, well-posedness }
\vspace{0.5cm}

\maketitle
	\renewcommand{\thefootnote}{} 
	\footnotetext{\noindent Department of Mathematics, Indian Institute of Technology Bhilai, 491002, Chhattisgarh, India\\  
		Emails: tapaswinisahu70@gmail.com, swati@iitbhilai.ac.in (Corresponding author)}
	\vspace{-2em}

\noindent\textbf{ABSTRACT.}
  In this paper, we investigate the numerical solution of a space-time fractional diffusion equation, involving Caputo fractional derivative $_{}^CD_t^\alpha$ in time and a Riesz fractional Laplacian $(-\Delta)^{\beta/2}$ in space. An $L1$ type discretization is employed for the temporal fractional derivative, together with a suitable finite difference spatial discretization of the fractional Laplacian. A graded temporal mesh is used to efficiently handle the weak singularity near the initial time. Wellposedness and error estimates are established in $L^2(\Omega)$ and fractional Sobolev norms. Numerical experiments are carried out for several test problems to examine the accuracy and convergence behavior of the proposed method. The numerical results support the theoretical estimates and demonstrate the effectiveness of the scheme.
\vspace{0.2cm}


\section{Introduction}

Nonlocal partial differential equations arise in science, engineering, and even finance, wherever nonlocal effects are modeled through nonlocal operators. Such problems were initially treated using integro-differential operators, as in the nonlocal Allen-Cahn equation \cite{du2019maximum}. Over time, however, non locality in \textit{time}, expressed through \textit{memory} effects, gained prominence, and this gave rise to a new branch of mathematics: fractional calculus.

In recent decades, the fractional derivatives have emerged an important area of research owing to their capability to model the hereditary properties and memory present in many biological \cite{magin2006fractional}, physical \cite{metzler2000random}, control-theoretic \cite{sabatier2007advances,naifar2022fractional}, and engineering systems \cite{sun2019review}, which are described more accurately by fractional models than by classical ones. Finite difference schemes are simple to implement and computationally efficient, making them attractive for obtaining accurate numerical solutions to a wide range of problems.

The development of fractional diffusion equations can be traced through a gradual evolution of fractional calculus from a primarily theoretical concept \cite{kilbas2006theory} to a well-established framework for modeling nonlocal and anomalous phenomena. Significant mathematical advances in the latter part of the 20th century established a rigorous foundation for fractional derivatives and their applications \cite{podlubny1998fractional}. Since then, considerable research has been devoted to the well-posedness, regularity, and stability of fractional models, as well as to the design and analysis of accurate and efficient numerical methods.

The time-fractional derivative introduces memory effects, that means the present state of the system depends upon its past behavior. In particular, Caputo fractional derivative of a function at any given time depends over the entire preceding time interval. This makes the process non-Markovian. In contrast, the space-fractional derivative captures long-range interactions, reflecting the non-local nature of the spatial diffusion. These properties make fractional models useful for describing anomalous-diffusion and other processes in which the systems behavior is influenced by its history and spatial interactions. 

Fractional-diffusion equations with time-, space-, or space-time-fractional derivatives arise naturally in the modeling of real-world transport phenomena. Time-fractional diffusion equations are commonly use to describe sub-diffusive process in which particle transport is hindered by trapping effects or long waiting times, as observed in disordered systems and heterogeneous materials \cite{metzler2000random,metzler2002space}. On the other hand, space-fractional diffusion equations describe supper diffusive transport characterized by long-range jumps, often modeled using L\'evy flights and have applications in transport through heterogeneous porous media \cite{metzler2000random,berkowitz2002physical}. Space-time-fractional diffusion equations provide a unified framework for modeling anomalous transport involving both mechanisms \cite{metzler2002space,metzler2004restaurant}.
The continuous-time random walk (CTRW) provides a microscopic description of particle diffusion in heterogeneous media, from which macroscopic models, such as the time-fractional reaction-diffusion equation, can be derived in the scaling limit \cite{metzler2000random,berkowitz2002physical}.

Recently, the study of time-space fractional diffusion equation has gained a significant attention. Mahmoud A Zaky in \cite{zaky2025finite}, has worked on the numerical solution of time-space fractional reaction diffusion equation with the help of Petrov-Galerkin spectral method. A discrete comparison principle \cite{chen2020discrete} by Martin Stynes and Hu Chen, fast evaluation of Caputo fractional derivative \cite{jiang2017fast}, an inverse problem \cite{zhang2018bayesian}, modeling and simulation \cite{gomez2016modeling}, like wise many work \cite{tan2022l1,gu2021fast,arshad2018finite,liu2026modified,hu2018finite}, has been done and numerous researchers are still contributing to the theoretical and computational field of diffusion equations. Despite of this, the existing well-posedness theory for Caputo-type fractional evolution equations have been defined predominantly on fixed spatial domains employing analytical tools, such as compactness method, fixed-point arguments, fractional energy estimates (see, \cite{diethelm2002analysis,diethelm2010analysis,van2021well,hamdi2026well,abdelaziz2023well}). The regularity of exact solution plays a fundamental role in determining the convergence rate and behavior of the numerical scheme. Researchers in \cite{duo2018novel,duo2019accurate}, have demonstrated the regularity and convergence behavior under various assumption on the initial data, source term, fractional operators domain, etc.

Although finite difference method have been widely investigated for space-time fractional diffusion equations, research has also progress towards finite element and spectral approach. For instance, \cite{feng2018unstructured} developed an unstructured mesh finite element method for two dimensional time-space Riesz fractional diffusion equation on irregular domains, establishing stability and convergence result. Similarly Li and Xu \cite{li2009space} proposed a space–time spectral method for time-fractional diffusion equations and established convergence estimates. These studies illustrate the continuing development of numerical techniques for fractional diffusion problems. Nevertheless, challenges associated with the nonlocal nature of fractional operators, solution regularity, computational efficiency, and rigorous convergence analysis continue to motivate further investigation.
However, to the best of our knowledge, no published work has yet addressed the well-posedness and the numerical approximation of the fractional diffusion equations, using finite difference framework that combines the $L1$ approximation for the temporal Caputo derivative discretization and a suitable finite difference approximation of Riesz fractional Laplacian $(-\Delta)^{\beta/2}$, the nonlocal generalization of classical Laplacian, in space. This identified gap in the literature provides the main motivation for our present work.

The $L1$ approximation is one of the simplest and most widely used discretization techniques for approximating the Caputo fractional derivative $^CD_t^\alpha$, which introduces memory effect to the system \cite{stynes2022survey}. Among the various types of boundary conditions, the homogeneous Dirichlet boundary condition has received considerable attention owing to its mathematical tractability and physical relevance. We establish the existence and uniqueness of weak solutions with the help of Sobolev spaces, and developed a fully discrete finite difference method for its numerical approximation.
 In particular, the spatial domain is discretized using uniform mesh, while graded mesh is employed in temporal domain. Finally, we established error estimate for proposed scheme. 
 
This paper is organized as follows. Section \ref{Sec.Model} introduces the model problem. In Section \ref{Sec.Wellpose}, we discuss the existence and uniqueness of the problem. Section \ref{Sec.Num.sol.} reviews some preliminary results and numerical solutions. Error analysis and corresponding numerical results are shown in Section \ref{Sec.Err} and Section \ref{Sec.Examples}, respectively. Finally, Section \ref{Sec.Conclude} concludes the paper with a discussion of main findings and motivation for the future work.

\section{\textbf{Model problem}}\label{Sec.Model}

In this section, we introduce the mathematical problem under consideration. The model problem is taken from the existing literature \cite{hu2018finite} . For detailed understanding of the formulation and the underlying modeling aspects, we refer the reader to the relevant literature.\\
Let $ T>0$ and $\Omega= (a,b)$ be a one dimensional spatial domain. We consider the following diffusion equation
\begin{align}
\label{Main_Eq}
\begin{cases}
{}^{C}D_t^\alpha u(x,t) + (-\Delta)^{\beta/2}u(x,t)
= f(x,t), 
& \quad (x,t)\in\Omega\times(0,T],\\
u(x,t)=0,
& \quad (x,t)\in(\mathbb{R}\setminus\Omega)\times(0,T],\\
u(x,0)=\phi(x),
& \quad x\in\Omega,
\end{cases}
\end{align}
 for $\alpha\in(0,1)$ and $\beta\in(0,2)$. $f(x,t)$ denotes the source (forcing) term of the fractional differential equation, which represents the external inputs acting on the system. Whereas $\phi(x,t)$ is the initial condition, which specifies the value of the solution at the initial time $t=0$. To establish the existence and uniqueness of a weak solution for this model problem, we assume $f\in L^2(0,T;L^2(\Omega))$ and $\phi\in L^2(\Omega)$.
 These assumptions provide the natural $L^2$-framework for the weak formulation of the problem and ensure that both the initial data and the source term are controlled in the energy estimates.
 Further details are discussed in Section \ref{Sec.Wellpose}. 
 
The derivative in the time direction $ ^CD_t^\alpha u$  is taken as a Caputo fractional derivative of order $\alpha \in (0,1),$ and is given by
\begin{equation}
    _{}^CD_t^\alpha u (x,t)=\frac{1}{\Gamma(1-\alpha)}\int_0^t (t-s)^{-\alpha} u^1(s) ds,
\end{equation}	
where $u$ is absolutely continuous on $[0,T]$, i.e., $u\in AC((0,T])$. 

Unlike the classical first-order derivative, the Caputo derivative is a nonlocal operator in time, since its value depends on the history of the function over the entire interval $[0,t]$ for $t \leq T$. Hence, it is particularly suitable for modeling processes exhibiting memory and hereditary effects.\\
For space direction, we have fractional Laplacian operator $(- \Delta)^{\beta/2}$ for $\beta \in (0,2),$ defined by
\begin{align}\label{Frac-laplace operator}
    (- \Delta)^{\beta/2} u(x) =c_{1,\beta} \ \text{P.V.}\int_\mathbb{R}\frac{u(x)-u(y)}{|x-y|^{1+\beta}}dy,
\end{align}
where the normalizing constant $c_{1,\beta}$ is given by 
\[c_{1,\beta}=\frac{2^{\beta-1}\beta \Gamma(\frac{\beta+1}{2})}{\sqrt{\pi}\Gamma(1-{\beta/2})}.\]
$\Gamma()$ represents the gamma function and  P.V. indicates Cauchy principal value.
This equation is a mathematical model to describe the non-uniform diffusion of matter in space-time. The Caputo derivative represents temporal memory effects, while the fractional Laplacian operator accounts for the non-local spatial diffusion. Therefore, the combination of these two operators provides a mathematical model for describing long-range interactions and anomalous diffusion involving both spatial and temporal non-locality. A diffusion equation incorporating both time and space fractional derivatives has been widely employed to anomalous diffusion phenomena. Similar diffusion problem have been investigated in the literature  by several authors like \cite{zaky2025finite,gu2021fast,yang2011novel} using different schemes and analysis. This inherent non-locality of fractional operators introduces significant challenges in both theoretical and numerical computations.

\section{\textbf{Well-posedness of the Model Problem}}\label{Sec.Wellpose}
We begin this section by recalling some standard spaces with usual norms and associated results. The  preliminary will be used while establishing the existence and uniqueness of the problem \eqref{Main_Eq}.

\begin{itemize}
    \item The space $L^2(\Omega)$ consisting of all square integrable measurable functions over $\Omega$ is defined by
\[ L^2(\Omega)=\left\{ w : \Omega \to \mathbb{R} : \int_\Omega |w(x)|^2 dx< \infty \right\},
\]
and is equipped with the norm
\[\norm{w}_{L^2(\Omega)}= \left( \int_\Omega |w(x)|^2dx\right)^{1/2}.
\]

\item The fractional Sobolev space $H^{\beta/2}(\Omega)$, of order $\beta/2$, for $\beta\in(0,2)$ is defined as 
\[
H^{\beta/2}(\Omega)=\left\{ w\in L^2(\Omega):\int\int_{\Omega\times\Omega}\frac{|w(x)-w(y)|^2}{|x-y|^{1+\beta}}dxdy<\infty\right\},
\] 
equipped with \textit{Slobodetskii norm} \cite{bhattacharyya2012distributions} 
\[
\norm{w}_{H^{\beta/2}}(\Omega)=\left[\norm{w}_2^2+\int\int_{\Omega\times\Omega}\frac{|w(x)-w(y)|^2}{|x-y|^{1+\beta}}dxdy\right]^{1/2}.
\]

\item The natural energy space $\widetilde{H}^{\beta/2}(\Omega)$ is defined as
\[
\widetilde{H}^{\beta/2}(\Omega)=\{w\in H^{\beta/2} : w=0 \quad \text{a.e.}\ \text{in}\ \mathbb{R}\setminus\Omega\}.
\]
\end{itemize}
For the additional details, one may follow \cite{yuan2026wellposedness,bhattacharyya2012distributions}.

For keeping the representation simple, we denote the above Hilbert spaces by the following notations: 
$$\mathcal{H}:=L^2(\Omega) ,\ V:=\widetilde{H}^{\beta/2}(\Omega),$$
with the corresponding usual norms. 
Let $V^*$ denotes the dual of V. The spaces $V,\ \mathcal{H},$ and $V^*$ form a Gelfand triple \cite{hamdi2026well} 
$$V\hookrightarrow \mathcal{H}\hookrightarrow V^*,$$
with continuous and dense embedding. 

We also represent the fractional Laplacian operator $(-\Delta)^{\beta/2}$ by $K$. For theoretical purpose that we see in the following sections, we restrict the domain of $K$ to $\widetilde{H}^\beta$ and denote this restricted domain by $D(K)$, and hence, we have \cite{yuan2026wellposedness}
$$K:D(K)\subset \mathcal{H} \to \mathcal{H}.$$

\subsection{Definition}(Weak Solution, \cite{hamdi2026well}).\label{Weak def}
We say $w$ is a weak solution of \eqref{Main_Eq} if 
\[w\in L^2(0,T;V),\quad ^CD_t^\alpha w \in L^2(0,T;V^*),\]and 
\[< {}^CD_t^\alpha w(t),v>_{V^*,V}+a(w(t),v)_\mathcal{H}=(f(.,t),v)_\mathcal{H},\quad\text{for all} \;v\in V,\] for almost every $ t\in(0,T),$ with the initial condition $w(.,0)=\phi(x).$

\subsection{Existence and Uniqueness.}\label{Existence}

Throughout this work, we impose the following assumptions, which enable us to use various inequalities to establish the existence and uniqueness of the considered problem:
\begin{enumerate} \label{ (A1)-(A3)}
\item[(A1)] $\phi(x)\in \mathcal{H};$

\item[(A2)] $f$ is measurable in $(x,t);$

\item[(A3)] $f\in L^2(0,T;\mathcal{H}).$ 
\end{enumerate}

As $K$ is a positive self adjoint operator with compact inverse, by spectral theory \cite{helffer2013spectral},  $K$ admits a sequence of orthonormal eigen pairs 
$\{(\lambda_ n,\psi_ n)\}_{n\geq1},$ such that, 
$ K \psi_n=\lambda_n \psi_n,$ $  
 0<\lambda_1 \leq \lambda_2 \leq \cdots \leq\lambda_n \to\infty.$ 
Now, consider our model problem
\[^CD_t^\alpha u(x,t) +Ku(x,t) =f(x,t),\quad x\in \Omega,\quad t\in(0,T],\]
by applying the Laplace transformation, the above equation becomes 
\begin{align}\label{lapl.1}
   (S^\alpha+\lambda_n)U_n(S)=S^{\alpha-1}\phi_n+F_n(S), 
\end{align}
where $U_n(S)$ and $F_n(S)$ represent the Laplace transforms of $u_n(t)$ and ${f_n}(t),$ respectively. Consequently, 
\begin{align}\label{Sum_e.f}
    u=\sum_{n=1}^\infty u_n\psi_n, \qquad
f=\sum_{n=1}^\infty f_n\psi_n, \qquad
\phi=\sum_{n=1}^\infty \phi_n\psi_n.
\end{align} 

\noindent After applying the inverse Laplace on Eq. \eqref{lapl.1}, we get
\[u_n=E_\alpha(-\lambda_nt^\alpha)\phi_n+\int_0^t (t-s)^{\alpha-1}E_{\alpha,\alpha}(-\lambda_n(t-s)^\alpha) f_n(s)ds.\]
Now by multiplying the above equation with $\psi_n$,  summing from $n=1$ to $\infty$, and substituting Eq. \eqref{Sum_e.f} give the following formulation 
\begin{align}\label{laplace sol}
u(t)=M_\alpha(t)\phi+\int_0^t N_\alpha(t-s)f(s)ds,\quad t\in[0,T],
\end{align}
where
\[
M_\alpha(t):=\sum_{n=1}^\infty E_\alpha(-\lambda_n t^\alpha)(.,\psi_n)\psi_n    ,\quad   N_\alpha(t):=\sum_{n=1}^\infty t^{\alpha-1}E_{\alpha,\alpha}(-\lambda_n t^\alpha)(.,\psi_n)\psi_n
.\]
Here, $E_\alpha$ and $ E_{\alpha,\alpha}$ represent the Mittag-Leffler functions of one parameter and two parameter, respectively. In general, the two parameter Mittag-Leffler function is defined as 
$$ E_{\alpha,\beta}(Z)=\sum_{k=0}^\infty\frac{Z^k}{\Gamma(\alpha k+\beta)},$$  
and for $\beta=1$ it reduces to $E_\alpha$.

\noindent Since $M_\alpha(t)$ and $ N_{\alpha}(t)$ form bounded operator families (see \cite{hamdi2026well}), hence $u\in C([0,T];L^2(\Omega))$.
\begin{theorem}(Existence).\label{Theorem-existence}
    Let assumptions $(A1)-(A3)$ of \eqref{Existence} hold true, then for every $\phi\in \mathcal{H},$ there exists a unique solution $u \in L
    ^2(0,T,V)\cap C([0,T];\mathcal{H})$, solving the formulation Eq. \eqref{laplace sol}.
\end{theorem}

 \noindent By the definition of weak solution \eqref{Weak def}, we see that the solution $u$ by Theorem \ref{Theorem-existence} forms a weak solution on $[0,T]$.

\begin{proof}
As  $$(- \Delta)^{\beta/2} u(x) =c_{1,\beta}\; \text{P.V.}\int_R\frac{u(x)-u(y)}{|x-y|^{1+\beta}}\;dy,$$
we start with considering the approach based on \cite{bonito2017numerical}, and for that, we assume
    \begin{align}\label{bilinear form}
    a(u,v)= \frac{c_{1,\beta}}{2}\int_\Omega\int_\Omega\frac{(u(x)-u(y))(v(x)-v(y))}{|x-y|^{1+\beta}}\;dxdy.
    \end{align}
One can easily show that $a(u,v)$ forms a bilinear form \cite{zhang2018riesz}. Applying the Cauchy-Schwarz inequality on Eq. \eqref{bilinear form}, makes the bilinear form $a(u,v)$ continuous on $V\times V,$ and hence, there exists a $M>0$ such that 
$$|a(u,v)|\leq M \norm{u}_V \cdot \norm{v}_V, \qquad \text{for all} \; u,v\in V.$$
To obtain the weak formulation, we first observe that $u\in L^2(0,T;V)$ implies that $Ku\in L^2(0,T;V^*). $ Consequently  $K:V\to V^*$ forms a bounded linear operator. With this using Riesz's theorem \cite{helffer2013spectral}, we have 
\begin{align*}
   <Ku,v>_{V^*,V} = a(u,v)_{L^2(\Omega)}, 
\end{align*}
therefore
\begin{align*}
    |<Ku,v>_{V^*,V}| =|a(u,v)| & \leq M\norm{u}_V\cdot\norm{v}_V,\quad\text{for all}\; u,v\in V, \\
    \implies \frac{|<Ku,v>_{V^*,V}|}{\norm{v}_V} & \leq M\norm{u}_V, \quad \text{holds}\quad  \text{for all}\; v\neq 0, \\
     \implies  \text{sup}_{v\neq 0} \frac{|<Ku,v>_{V^*,V}|}{\norm{v}_V} & \leq M \norm u_V .
\end{align*}
The above estimate yields the inequality
\begin{align}\label{lap_norm}
 \norm {Ku}_{V^*}\leq M\norm{u}_V .
\end{align}
It follows that for $u\in L^2(0,T;V)\cap C([0,T];\mathcal{H})$, $Ku\in L^2(0,T,V^*)$. Now, since $\mathcal{H}\hookrightarrow V^*$ is a continuous embedding, from (A3) we obtain
\begin{align}
    f\in L^2(0,T;\mathcal{H})\subset L^2(0,T;V^*) .\label{forcing norm}
\end{align}
Together with the Eq. \eqref{lap_norm} and Eq. \eqref{forcing norm}, the equation $^C D_t^\alpha u=f-Ku$ is satisfied in $L^2(0,T;V^*).$
Multiplying with an arbitrary test function $v\in V$, $u$ satisfies the weak formulation
\begin{align} \label{weak formulation}
<{}^CD_t^\alpha u,v>_{V^*,V}+a(u,v)_\mathcal{H}=(f,v)_\mathcal{H} .
\end{align}
The above result shows the existence of weak solution of \eqref{Main_Eq}. 
\end{proof}

We now discuss the uniqueness of the weak solution in the form of following theorem.

\begin{theorem}(Uniqueness).
    Assume $(A1)-(A3)$ of \eqref{Existence} hold. Let $w_1$ and $w_2$ denote weak solutions of Eq. \eqref{Main_Eq} on [0,T] associated with the initial data $\phi_1,\phi_2\in \mathcal{H}$, then\[
    \norm{w_1(t)-w_2(t)}_\mathcal{H}^2\leq\norm{\phi_1-\phi_2}_\mathcal{H}^2, \quad \text{for all}\; t\in[0,T].
    \]
    In particular, if $\phi_1=\phi_2,$ implies $w_1=w_2$ on [0,T].
\end{theorem}

\begin{proof}
Being weak solutions of \eqref{Main_Eq}, $w_1$ and $w_2$ satisfy the weak formulation \eqref{weak formulation}, that is, 
\[
<{}^C D_t^\alpha w_1,v>_{V^*,V}+a(w_1,v)_\mathcal{H}=(f,v)_\mathcal{H},\]
\[<{}^C D_t^\alpha w_2,v>_{V^*,V}+a(w_2,v)_\mathcal{H}= (f,v)_\mathcal{H}.
\]
Subtracting the above formulations, replacing  $w_1-w_2$ by $u$, and also testing with $v=u\in V$, yield
\[ <{}^CD_t^\alpha u(t),u(t)>_{V^*,V}+\norm{K^{1/2}u(t)}_\mathcal{H}^2= 0.\]

\noindent Using fractional energy identity Lemma 4.2 of \cite{hamdi2026well}, we obtain
\[\frac{1}{2}\;{^CD_t^\alpha }\norm{u(t)}_\mathcal{H}^2+\norm{K^{1/2}u(t)}_\mathcal{H}^2\leq0 .\]
By dropping $\norm{K^{1/2}u(t)}_\mathcal{H}^2$ for its non-negative behavior, and setting $y(t):=\norm{u(t)}_\mathcal{H}^2$, we obtain
\[\frac{1}{2}\; {^CD_t^\alpha} y(t)\leq0.\]
Further by similar approximation to Lemma 5.4 \cite{hamdi2026well} implies 
$$y(t)\leq y(0),\quad  y(0)=\norm{\phi_1-\phi_2}_\mathcal{H}^2 ,\quad t\in[0,T].$$
Hence, 
\[ \norm {w_1(t)-w_2(t)}_\mathcal{H}^2\leq\norm{\phi_1-\phi_2}_\mathcal{H}^2.\]
In particular, if $\phi_1=\phi_2,$ the required estimate follows.
\end{proof}
\section{Numerical Discretizations for the Model problem}\label{Sec.Num.sol.}

The solution of the model problem \eqref{Main_Eq} is numerically approximated by using a finite difference scheme on a graded temporal mesh and a uniform spatial mesh. With the given positive integers $M$ and $N$, we define the spatial step size by $h=(b-a)/M$ and the spatial nodes by $x_i=ih,\; 0\leq i\leq M$. The temporal graded mesh is defined as
\[t_n=T\left(\frac{n}{N}\right)^r,\quad n=0,1 , \cdots, N,\]
where $r>1$ is the mesh grading parameter and for $r=1$, it becomes uniform mesh. Set the corresponding time-step size 
\[\tau_n=t_n-t_{n-1},\quad n=1,2 , \cdots, N.\]
Consequently, the computational mesh is given by 
$\{(x_i,t_n):i=0,1, \cdots, M ,\;n=0,1, \cdots, N\}\;.$
\subsection{Numerical approximation of Caputo derivative:} 

The Caputo fractional derivative $^CD_t^\alpha$ is approximated by $L1$ approximation \cite{stynes2017error} as follows 
\begin{align}\label{L1 discretization}
    {^CD_N^\alpha} u_i^n & =\frac{1}{\Gamma(1-\alpha)}\sum_{j=0}^{n-1} \frac{u_i^{j+1}-u_i^j}{\tau_{j+1}}\int_{s=t_j}^{t_{j+1}}(t_n-s)^{-\alpha}ds, \nonumber \\
   &= \frac{1}{\Gamma(1-\alpha)}\sum_{j=0}^{n-1}\frac{u_i^{j+1}-u_i^j}{\tau_{j+1}}\left[\frac{(t_n-s)^{1-\alpha}}{-(1-\alpha)}\right]_{t_k}^{t_{k+1}},\nonumber \\
  & =\frac{1}{\Gamma(2-\alpha)}\sum_{j=0}^{n-1}\frac{u_i^{j+1}-u_i^j}{\tau_{j+1}}[(t_n-t_j)^{1-\alpha}-(t_n-t_{j+1})^{1-\alpha}],\nonumber\\
 & =\frac{b_{n,1}}{\Gamma(2-\alpha)}u_i^n-\frac{b_{n,n}}{\Gamma(2-\alpha)}u_i^0+\frac{1}{\Gamma(2-\alpha)}\sum_{j=1}^{n-1}u_i^{n-j}[b_{n,j+1}-b_{n,j}],
\end{align}

\noindent where
\[b_{n,j}: =\frac{(t_n-t_{n-j})^{1-\alpha}-(t_n-t_{n-j+1})^{1-\alpha}}{\tau_{n-j+1}},\quad \text{for}\;j=1,2 , \cdots, n-1.\]
For $j=1$, $b_{n,1}=\tau_n^{-\alpha}.$\\

\subsection{Temporal error estimate:} \label{temp-err}

  For time fractional diffusion problem the solution generally exhibits a weak singularity near the initial time at $t=0$. 
On the uniform mesh, this limited regularity generally leads to the convergence rate of $O(N^{-\alpha}).$ To improve the temporal accuracy, graded mesh is employed, that provides a finer time-step distribution near the initial point.
Stynes, O'Riordan and Gracia in \cite{stynes2017error}, developed an error analysis of $L1$ finite difference approximation for $i\in \{1,2, \cdots, M-1\}$ and $n\in \{1,2, \cdots, N\}$, which states that there exists a constant $ C$ such that
\[
|D_N^\alpha u_i^n-D_t^\alpha u_i^n|\leq C N^{-min\{r\alpha,\ 2-\alpha\}},\quad \text{for all}\; (x_i,t_n)\in (0,l)\times(0,T]\;.
\]
The authors employed both uniform and graded mesh that establish consistency and convergence estimates for $L1$ scheme.
Setting $r=(2-\alpha)/\alpha,$ the method achieves optimal convergence rate $O(N^{-(2-\alpha)}).$ This estimate indicates that the temporal accuracy depends on both the fractional order $\alpha$ and the choice of graded mesh $r.$

\subsection{Numerical approximation of Laplacian operator.} 
The approximation of the fractional Laplacian operator $(-\Delta)^{\beta/2}u_i^n$ by \cite{xing2023numerical} is given as
\begin{align}\label{spatial-app}
(-\Delta)^{\beta/2}u_i^n \approx \delta_{h,\mu}^\beta u_i^n  = c_{}^{\beta,\mu} \sum_{v=1}^{M-1} g_{i-v}^{(\beta)} u_v^n =AU^n=(-\Delta)_{h,\mu}^{\beta/2}u, 
\end{align}
where
\[
c_{}^{\beta,\mu}=\frac{c_{1,\beta}}{p h^\beta},\ p= {\mu-\beta},\ \mu=1+{\beta/2}.
\]
The matrix $A$ is a symmetric, positive-definite and Toeplitz, with entries given by the product of $c_{}^{\beta,\mu} $ and the corresponding sum involving  $g_v^{(\beta)}$, where $g_v^{(\beta)}$ is defined as follows
\[g_v^{(\beta)}=
\begin{cases}
\sum_{l=2}^{M-1}\frac{(l+1)^p-(l-1)^p}{l^\mu}+\frac{M^p-(M-1)^p}{M^\mu}+2^{p}+\frac{2 p}{\beta M^\beta}, & v= 0,\\
-2^{p-1}, & v=1,\\
-\frac{(v+1)^p-(v-1)^p}{2v^\mu}, &v=2,\cdots, M-2,\\
g_{-v}^{(\beta)},&v=-1,-2,\cdots, 2-M.
\end{cases}
\]

 \subsection{Spatial error estimate.}\label{spatial err}
 To examine the spatial convergence of the finite difference approximation  $(-\Delta)_{h,\mu}^{\beta/2}$ 
 to the fractional Laplacian $(-\Delta)^{\beta/2}$, several cases are considered based on the regularity of the solution and the choice of \textit{splitting parameter} $\mu$. In particular convergence behavior depends on both fractional order $\beta$ and the value of $\mu$. For the deeper understanding of the statements made, the reader may follow  \cite{duo2018novel}.\\
\textbf{R1:} For $u\in C^{1,\beta/2}(\mathbb{R}),$ and for $\mu\in(\beta,2],$ there is a positive constant $C$ depending on $\beta$ and $\mu$ such that, the approximation satisfies
\[ 
\norm{(-\Delta)^{\beta/2}u-(-\Delta)_{h,\mu}^{\beta/2}u}_\infty\leq Ch^{1-\beta/2},\quad \beta\in(0,2).
\]
\textbf{R2:} For $u\in C^{3,\beta/2}(\mathbb{R}),$ with finite support on $\Omega\subset\mathbb{R},$ there is a positive constant $C$ depending on $\beta\;\text{and}\; \mu$ such that the approximation attains second order accuracy
\[
\norm{{(-\Delta)^{\beta/2}u-(-\Delta)_{h,\mu}^{\beta/2}u}}_{\infty,\Omega}\leq Ch^2,\quad \beta\in(0,2) , 
\]
\noindent for the splitting parameter $\mu=2$ or $1+{\beta/2}.$\\
\textbf{R3:}
For $u\in C^{3,\beta/2}(\mathbb{R}),$ when $\mu \in(\beta,2)$ with $\mu\neq \{1+\frac{\beta}{2} ,2\}$ the approximation has convergence order $O(h^{2-\beta})$ and the following inequality holds
\[
\norm{{(-\Delta)^{\beta/2}u-(-\Delta)_{h,\mu}^{\beta/2}u}}_{\infty,\Omega}\leq Ch^{2-\beta},\quad \beta\in(1,2) .
\]
For $ \beta\in(0,1]$, establishing uniform error estimate valid for all $\mu\in(\beta,2) $ is considerably more involved.
Therefore an appropriate selection of $\mu$, together with the fractional order $\beta/2$, provides useful information about the expected error and the convergence rate when presenting and interpreting the numerical results.\\[2mm]
\textbf{Remark:} The constant $C$ appearing above denotes a generic positive real constant, which may take different values at different occurrences.
\subsection{Numerical solution of model problem.}

Substituting Eq. \eqref{L1 discretization} and Eq. \eqref{spatial-app} in Eq. \eqref{Main_Eq} we obtain,
\[ 
c_{}^{\beta,s}\sum_{k=1}^{M-1}g_{i-k}^{(\beta)} u_k^n+ \frac{b_{n,1}}{\Gamma(2-\alpha)}u_i^n-\frac{b_{n,n}}{\Gamma(2-\alpha)}u_i^0+\frac{1}{\Gamma(2-\alpha)}\sum_{j=1}^{n-1}u_i^{n-j}[b_{n,j+1}-b_{n,j}]=f_i^n.
\]
By rearranging terms, we have
\[
c_{}^{\beta,s}\sum_{k=1}^{M-1}g_{i-k}^{(\beta)} u_k^n+ \frac{b_{n,1}}{\Gamma(2-\alpha)}u_i^n=f_i^n+\frac{b_{n,n}}{\Gamma(2-\alpha)}u_i^0-\frac{1}{\Gamma(2-\alpha)}\sum_{j=1}^{n-1}u_i^{n-j}[b_{n,j+1}-b_{n,j}].
\]
Rewriting the above equation into matrix form, we obtain the numerical scheme to solve the problem \eqref{Main_Eq}
\begin{align}\label{n_sol}
 G_1U^n &=F^n+\frac{b_{n,n}}{\Gamma(2-\alpha)}U^0-\frac{1}{\Gamma(2-\alpha)}\sum_{j=1}^{n-1}U^{n-j}[b_{n,j+1}-b_{n,j}],
\end{align}
where
\[
G_1=c^{\beta,s}
\begin{bmatrix}
    l_n & g_{-1}^{(\beta)} & g_{-2}^{(\beta)}&  \cdots &g_{-M+2}^{(\beta)}\\
    g_1^{(\beta)} &l_n & g_{-1}^{(\beta)}&  \cdots &g_{-M+3}^{(\beta)}\\
    \vdots& \ddots & \ddots &\ddots & \vdots \\
     g_{M-3}^{(\beta)} & \ddots & \ddots  &l_n & g_{-1}^{(\beta)}\\
    g_{M-2}^{(\beta)} & g_{M-3}^{(\beta)} & \cdots & g_1^{(\beta)} &l_n
\end{bmatrix}_{M-1\times M-1},
\]
and  $l_n=g_0^{\beta}+\frac{\tau_n^{-\alpha}}{c^{\beta,s}\Gamma(2-\alpha)},\quad U^n=[u_1^n,u_2^n,\cdots,u_{M-1}^n]^T,\quad F^n=[f_1^n,f_2^n,\cdots,f_{M-1}^n]^T .$\\[2mm]
Since $G_1$ is a non singular symmetric Toeplitz matrix Eq. \eqref{n_sol} yields,

\begin{align}\label{solution}
    U^n=G_1^{-1}[F^n+\frac{b_{n,n}}{\Gamma(2-\alpha)}\phi(x)-\frac{1}{\Gamma(2-\alpha)}\sum_{j=1}^{n-1}U^{n-j}(b_{n,j+1}-b_{n,j})],
\end{align}
is the required solution for the problem \eqref{Main_Eq}.

\section{\textbf{Error Analysis.}}\label{Sec.Err}
\begin{theorem}
Let $u(x,t)$ be the solution of the finite difference scheme of the problem \eqref{Main_Eq} and $\{u_m^n\;|0\leq m\leq M, 0\leq n\leq N\} $ be the solution of the numerical scheme \eqref{solution} on the temporal mesh .Then there exists a positive constant C such that the total truncation error of the scheme is
\begin{align}
\|e^n\| \le C \left( h^q + N^{-\min\{r\alpha,\, 2-\alpha\}} \right),\;\; q\in(0,2],\;\; n = 1,2,\dots,N,
\end{align}
where
$e_m^n = U_m^n - u_m^n,\;\; 0 \le m \le M,\;\; 0 \le n \le N$.
\end{theorem}

\begin{proof}
The local truncation error quantifies the discrepancy between the exact solution and the numerical approximation at each grid point. For our numerical scheme \eqref{L1 discretization} and  \eqref{spatial-app},  we have
$|D_N^\alpha u_m^n- {}^CD_t^\alpha u_m^n|\leq CN^{-min\{2-\alpha,\;r\alpha\}}$ \cite{stynes2017error} and  $|(-\Delta)^{\beta/2}u_m^n-\delta_{h,\mu}^\beta u_m^n|=O(h^q),\;q\in(0,2]$ \cite{xing2023numerical}.
Define the error as
\begin{align*}
    e_m^n = U_m^n - u_m^n, 0\leq m\leq M ,0\leq n\leq N .
\end{align*}
Considering the Eq. \eqref{Main_Eq} at $(x_m,t_n)$, we obtain
\begin{align}\label{trunc_err}
    ^CD_N^\alpha U_m^n+\delta_{h,\mu}^{\beta} U_m^n = f_m^n+R_m^n,
\end{align}
\noindent where 
\[
R_m^n\leq C(h^q+N^{-min\{2-\alpha,\;r\alpha\}}).
\]
Now replacing numerical solution $u_m^n$ to replace $U_m^n$ and omitting $R_m^n$ in Eq. \eqref{trunc_err} we get 
\begin{align}
\begin{cases}
    ^CD_N^\alpha u_m^n+\delta_{h,\mu}^{\beta} u_m^n = f_m^n,\\
    u_0^n=0,\; u_M^n=0,\; 0\leq n\leq N,\\
    u_m^0 =\phi(x_m),\; 1\leq m\leq M-1.
\end{cases}
\end{align}
One can easily obtain the following error system of equations:
\begin{align}\label{main_err}
\begin{cases}
^CD_N^\alpha e_m^n + \delta_{h,\mu}^{\beta} e_m^n = R_m^n,\\
e_0^n = 0, \quad e_M^n = 0, \quad 0 \le n \le N,\\
e_m^0 = 0, \quad 1 \le m \le M-1.
\end{cases}
\end{align}
 Let $e^n = (e_1^n, e_2^n, \dots, e_{M-1}^n)^T$. Then the error equation in vector form is 
\begin{align}\label{error_1}
^CD_N^\alpha e^n + \delta_{h,\mu}^{\beta} e^n = R^n. 
\end{align}
Multiplying Eq. \eqref{main_err} by $h e_m^n$ and summing over all spatial nodes $m = 1 \;\text{to}\; M-1,$ we obtain
\begin{align}\label{error_2}
    (^CD_N^\alpha e^n, e^n) + (\delta_{h,\mu}^{\beta} e^n, e^n) = (R^n, e^n).
\end{align}
\noindent Firstly, the discrete Caputo derivative ($L1$ scheme) is given by (see \cite{liao2018sharp})
\[
    ^CD_N^\alpha e^n = \sum_{k=1}^{n} a_{n-k}^{(n)} (e^k - e^{k-1}),
\]
where $a_{n-k}^{(n)} > 0,$ and $a_{n-k-1}^{(n)} > a_{n-k}^{(n)}.$
Now
\begin{align}\label{error_3}
(^CD_N^\alpha e^n, e^n)
&= \sum_{k=1}^{n} a_{n-k}^{(n)} \left( e^k - e^{k-1}, \, e^n \right)\nonumber\\
&\ge \frac{1}{2} {^CD_N^\alpha} \|e^n\|^2.
\end{align} 
By direct computation, we obtain $(\delta_{h,\mu}^{\beta} e^n, e^n) = h (e^n)^T A e^n>0,$
where $A$ is symmetric positive definite matrix we have.\\
Now estimate of the right-hand side of Eq. \eqref{error_2}
$$(R^n, e^n) = h \sum_{m=1}^{M-1} R_m^n e_m^n.
$$\\
Applying the Cauchy--Schwarz inequality, we obtain
\begin{align}\label{error_4}
     |(R^n, e^n)| \le \|R^n\| \, \|e^n\|.
\end{align}
By combining all the above bounds of Eq. \eqref{error_2}, Eq. \eqref{error_3} and Eq. \eqref{error_4} we get
\begin{align*}
    \frac{1}{2} {^CD_N^\alpha} \|e^n\|^2 
\le \|R^n\| \, \|e^n\|,
\end{align*}
\[
\therefore
    {^CD_N^\alpha} \|e^n\|^2 
\le 2 \|R^n\| \, \|e^n\|.
\]\\
By using Young's inequality, we have
\begin{align*}
    {^CD_N^\alpha} \norm{e^n}^2 & \le  \left(  \norm{e^n}^2 + \|R^n\|^2 \right),
\end{align*}
Applying the discrete Gr\"onwall inequality \cite{liao2018sharp} for some positive constant $C_1$,  we obtain
\begin{align*}
     \norm{e^n}^2 &\le C_1\; \max_{1 \le j \le n} \|R^j\|^2,\\
    \|e^n\| =  &\le C_1\; \max_{1 \le j \le n} \|R^j\|.
\end{align*}
We know that the truncation error satisfies
\begin{align*}
    \|R^j\| \le C \left( h^q + N^{-\min\{2-\alpha,\; r\alpha\}} \right), 
\quad \text{for all}\; j.
\end{align*}
Hence
\begin{align*}
\|e^n\| \le C \left( h^q + N^{-\min\{2-\alpha,\; r\alpha\}} \right),\;n=1,2,\cdots,N.
\end{align*}
\end{proof}

\section{Numerical results}\label{Sec.Examples}

In this section, we present numerical results to validate the correctness and accuracy of the theoretical claims. In the first two examples, we are considering the domain $(-l,l)$, which is a standard choice of such kind of problems in the literature. In the subsequent two examples, we consider the domain $(0,b)$, that can be extended to a more general interval $(a,b)$.
In the initial three problems, where the exact solution is known, the error is calculated by comparing the numerical solution with the exact solution for different mesh sizes. The convergence rate is then estimated by comparing the errors obtained on successive meshes.\\
However, in the absence of exact solution, we compare the numerical solution obtained on two mesh sizes, usually a coarse mesh and a finer mesh. This approach is known as two-mesh principle.
Denote $u_i^n$ with $0\leq i\leq M,0\leq n\leq N $, the numerical scheme generated by Equation \eqref{solution} on the mesh introduced at the beginning of Section 4. Let $w_{2i}^{2n}$ denotes the numerical solution on the refined mesh with $2M$ spatial subdivisions and $2N$ temporal subdivisions, evaluated at the same physical grid point as of the coarser mesh.
The maximum difference between the two numerical solutions on these two meshes is defined as  
 \[
D_{M,N} := \max\limits_{\substack{0\leq i\leq M}}\max\limits_{\substack{0\leq n\leq N}}|w_{2i}^{2n}-u_i^n|.
\]
Similarly, $D_{2M,2N}$ denotes the maximum difference between the numerical solutions on the meshes $(2M,2N)$ and $(4M,4N)$.
The rate of convergence is then estimated as
\[
rate_{M,N}:= log_2\left(\frac{D_{M,N}}{D_{2M,2N}}\right).
\]
This approach allows the convergence rate to be estimated without explicitly knowing the exact solution.

\begin{example}\label{Example 6.1}
We consider a time-independent version of our model problem  \eqref{Main_Eq}, which reduces to a Poisson equation \cite{dyda2012fractional}, subject to a homogeneous Dirichlet boundary conditions:
\[
\begin{cases}
    (- \Delta)^{\beta/2} u(x)  = f(x) ,&x\in(-1,1),\\
u(x)=0,&x\in\mathbb{R}\setminus(-1,1).
\end{cases}
\]
For a particular forcing term $f(x)$ as
\[
f(x) = \frac{\Gamma(\beta+3)(\beta+4)}{120}\left(15-(10\beta+30)x^2+(\beta+3)(\beta+5)x^4\right)x,
\]
we get the exact solution  
 \[
  u(x)=\begin{cases}
     x(1-x^2)^{2+\beta/2},& \text{for}\; x\in \Omega=(-1,1),\\
      0, & \text{ for} \;x\in \mathbb{R}\setminus(-1,1).
 \end{cases}
 \]
The motivation behind choosing this example is to validate the rate  of convergence in space direction of the proposed numerical scheme.
    \begin{table}[htbp]
{ \scriptsize
\centering
\caption{Maximum errors and spatial convergence rate for Example \ref{Example 6.1} with $\mu=1+\beta/2$.}
\label{Table_1}
\begin{tabular}{c|cc|cc|cc|cc|cc}
\toprule
$\beta$
& \multicolumn{2}{c|}{$0.6$}
& \multicolumn{2}{c|}{$0.8$}
& \multicolumn{2}{c|}{$1$}
& \multicolumn{2}{c|}{$1.2$}
& \multicolumn{2}{c}{$1.6$} \\
\cmidrule(lr){2-3}
\cmidrule(lr){4-5}
\cmidrule(lr){6-7}
\cmidrule(lr){8-9}
\cmidrule(lr){10-11}

$M$
& Err & Rate
& Err & Rate
& Err & Rate
& Err & Rate
& Err & Rate \\
\midrule

32
& $2.3841e$-$04$ & --
& $2.9189e$-$04$ & --
& $3.4623e$-$04$ & --
& $4.0892e$-$04$ & --
& $6.3108e$-$04$ & -- \\

64
& $6.1106e$-$05$ & 1.96
& $7.0202e$-$05$ & 2.06
& $8.1309e$-$05$ & 2.09
& $9.3682e$-$05$ & 2.13
& $1.0460e$-$04$ & 2.17 \\

128
& $1.6106e$-$05$ & 1.92
& $1.7836e$-$05$ & 1.98
& $1.9727e$-$05$ & 2.04
& $2.2325e$-$05$ & 2.07
& $3.2251e$-$05$ & 2.12 \\

256
& $4.1380e$-$06$ & 1.96
& $4.5226e$-$06$ & 1.98
& $4.8621e$-$06$ & 2.02
& $5.4389e$-$06$ & 2.04
& $7.5477e$-$06$ & 2.10 \\

512
& $1.0485e$-$06$ & 1.98
& $1.1378e$-$06$ & 1.99
& $1.2072e$-$06$ & 2.00
& $1.3411e$-$06$ & 2.02
& $1.7931e$-$06$ & 2.07 \\
1024
& $2.6387e$-$07$ & 1.99
& $2.8526e$-$07$ & 2.00
& $3.0086e$-$07$ & 2.00
& $3.3280e$-$07$ & 2.01
& $4.3106e$-$07$ & 2.05 \\
\bottomrule
\end{tabular}
}
\end{table}
As clearly shown in Table \ref{Table_1}, the spatial convergence rate is $O(h^2)$, which is consistent with the theoretical result \eqref{spatial err} and demonstrate the optimal spatial convergence rate for different choice of $\alpha$ and $\beta$. Furthermore, we plot the corresponding numerical and exact solution in Figure \ref{figure1}, which shows the accuracy of the proposed method. We also present the spatial convergence rate graph in Figure \ref{figure1}, that demonstrates that the proposed method achieves second-order convergence rate, which is consistent with the theoretical prediction.

\begin{figure}[htbp]
    \centering
\label{Figure_1}
    \includegraphics[width=0.45\textwidth]{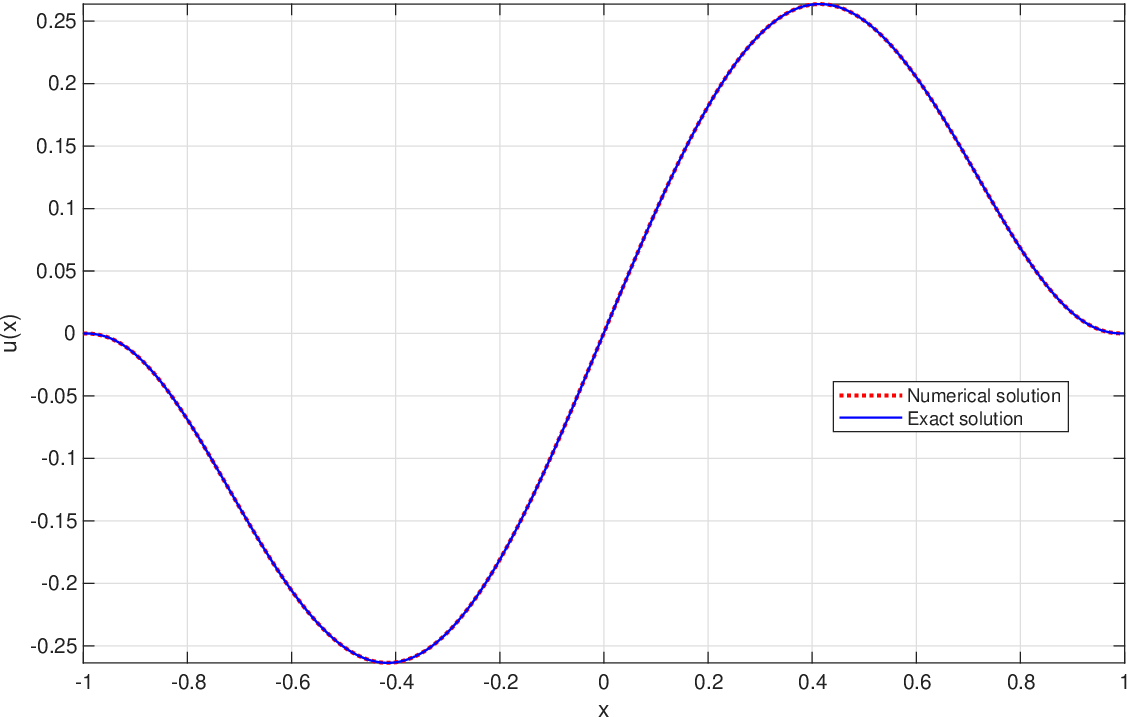}
    \hfill
    \includegraphics[width=0.45\textwidth]{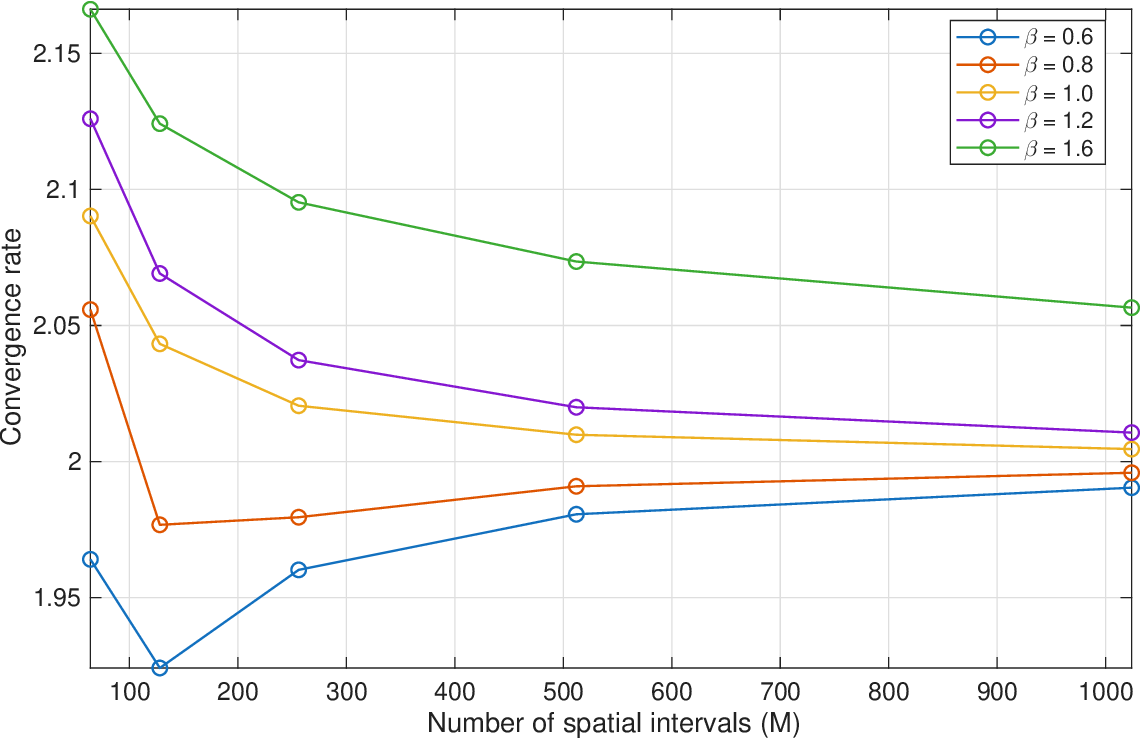}

\caption{The left above figure represents comparison between exact and numerical solutions for $\beta=0.8$ and $\mu=1+\beta/2,$ demonstrating the accuracy of the proposed numerical method, while the right figure illustrates the spatial convergence rate for different values of $\beta$ for Example  \ref{Example 6.1}.}
\label{figure1}
\end{figure}
\end{example}
\begin{example}\label{Example 6.2}
We consider a time-space fractional diffusion equation with prescribed exact solution, where the corresponding right hand side source term is derived from the governing equation by substituting the exact solution:
     \[
 \begin{cases}
      _{}^CD_t^\alpha u(x,t) + (- \Delta)^{\beta/2} u(x,t)  = f(x,t), \quad (x,t) \in (-1,1) \times (0,1],\\
      u(x,t)=0, \quad  x \in \mathbb{R}\setminus(-1,1),\\
      u(x,0)  = 0 , \quad  x \in (-1,1).
 \end{cases}
 \]
We choose the function then the exact solution for Example \ref{Example 6.2} is given as 
\[
u(x,t)  = t^2x(1-x^2)^{2+\beta/2}, \quad x\in(-1,1).
\]
The corresponding source function is given as
$f(x,t)$ for $x\in(-1,1)$ as $$f(x,t)= \frac{2t^{2-\alpha}}{\Gamma(3-\alpha)}x(1-x^2)^{2+\beta/2}+t^2\frac{\Gamma(\beta+3)(\beta+4)}{120}\left(15-(10\beta+30)x^2+(\beta+3)(\beta+5)x^4\right)x.$$

From Table\ \ref{Table_2}, we can clearly see that, for different choices of $\alpha$ and $\beta$, the graded time mesh with $r=\frac{2-\alpha}{\alpha}$ achieves the optimal convergence rate of $O(N^{-(2-\alpha)}),$ which agrees the theoretical error estimates \eqref{temp-err}. From Table \ref{Table_3}, the numerical result for the spatial case demonstrate that the proposed scheme achieve the optimal spatial convergence rate of $O(h^2)$ for different choices of fractional temporal and spatial orders. This confirms the second order accuracy of spatial discretization \eqref{spatial err} and shows that numerical results are in good agreement with the theoretical error estimate.
 \begin{table}[htbp]
{ \scriptsize
\centering
\caption{Maximum errors and temporal convergence rate for Example \ref{Example 6.2} with $\mu=1+\beta/2$,\; $r=(2-\alpha)/\alpha\;, \;\text{and} \; M=1024$.}
\label{Table_2}
\begin{tabular}{c|cc|cc|cc|cc|cc}
\toprule
$(\alpha,\beta)$
& \multicolumn{2}{c|}{$(0.3,0.5)$}
& \multicolumn{2}{c|}{$(0.5,1.2)$}
& \multicolumn{2}{c|}{$(0.7,1.4)$}
& \multicolumn{2}{c|}{$(0.8,1.6)$}
& \multicolumn{2}{c}{$(0.9,1.9)$} \\
\cmidrule(lr){2-3}
\cmidrule(lr){4-5}
\cmidrule(lr){6-7}
\cmidrule(lr){8-9}
\cmidrule(lr){10-11}

$N$
& Err & Rate
& Err & Rate
& Err & Rate
& Err & Rate
& Err & Rate \\
\midrule

32
& $9.3108e$-$04$ & --
& $5.4877e$-$04$ & --
& $6.2217e$-$04$ & --
& $6.2380e$-$04$ & --
& $5.5342e$-$04$ & -- \\

64
& $3.2385e$-$04$ & 1.52
& $2.0456e$-$04$ & 1.42
& $2.5753e$-$04$ & 1.27
& $2.7425e$-$04$ & 1.18
& $2.5940e$-$04$ & 1.09 \\

128
& $1.0880e$-$04$ & 1.57
& $7.5008e$-$05$ & 1.45
& $1.0589e$-$04$ & 1.28
& $1.2017e$-$04$ & 1.19
& $1.2157e$-$04$ & 1.09 \\

256
& $3.5780e$-$05$ & 1.60
& $2.7296e$-$05$ & 1.46
& $4.3443e$-$05$ & 1.29
& $5.2639e$-$05$ & 1.19
& $5.7107e$-$05$ & 1.09 \\

512
& $1.1639e$-$05$ & 1.62
& $9.9623e$-$06$ & 1.46
& $1.7873e$-$05$ & 1.29
& $2.3132e$-$05$ & 1.19
& $2.6997e$-$05$ & 1.08 \\
\bottomrule
\end{tabular}
}
\end{table}

\begin{table}[htbp]
{\scriptsize
\centering
\caption{Maximum errors and spatial convergence rate for Example \ref{Example 6.2} with $\mu=1+\beta/2$, $r=(2-\alpha)/\alpha$, and $N=2048$.}
\label{Table_3}
\begin{tabular}{c|cc|cc|cc|cc|cc}
\toprule
$(\alpha,\beta)$
& \multicolumn{2}{c|}{$(0.2,0.5)$}
& \multicolumn{2}{c|}{$(0.2,1)$}
& \multicolumn{2}{c|}{$(0.4,1.5)$}
& \multicolumn{2}{c|}{$(0.5,1.8)$}
& \multicolumn{2}{c}{$(0.6,1.9)$} \\
\cmidrule(lr){2-3}
\cmidrule(lr){4-5}
\cmidrule(lr){6-7}
\cmidrule(lr){8-9}
\cmidrule(lr){10-11}

$M$
& Err & Rate
& Err & Rate
& Err & Rate
& Err & Rate
& Err & Rate \\
\midrule

16
& $5.3351e$-$04$ & --
& $1.1845e$-$03$ & --
& $2.2303e$-$03$ & --
& $3.4776e$-$03$ & --
& $4.0749e$-$03$ & -- \\

32
& $1.5593e$-$04$ & 1.77
& $2.7562e$-$04$ & 2.10
& $4.7349e$-$04$ & 2.24
& $7.7389e$-$04$ & 2.17
& $9.4264e$-$04$ & 2.11 \\

64
& $4.2306e$-$05$ & 1.88
& $7.0597e$-$05$ & 1.97
& $1.1295e$-$04$ & 2.07
& $1.7988e$-$04$ & 2.11
& $2.2656e$-$04$ & 2.06 \\

128
& $1.1315e$-$05$ & 1.90
& $1.7657e$-$05$ & 2.00
& $2.6322e$-$05$ & 2.10
& $4.2335e$-$05$ & 2.09
& $5.4861e$-$05$ & 2.05 \\

256
& $3.0080e$-$06$ & 1.91
& $4.3995e$-$05$ & 2.00
& $6.3400e$-$06$ & 2.05
& $1.0387e$-$05$ & 2.03
& $1.3818e$-$05$ & 1.99 \\
\bottomrule
\end{tabular}
}
\end{table}
 Additionally, we have also plotted the corresponding graphs for selected parameter values in Figure \ref{figure2} to visualize the performance of the numerical scheme in both time and space. The spatial and temporal convergence plots provide further insight into the error reduction under mesh refinement and confirm the optimal convergence behavior of the numerical method in both dimensions.
\begin{figure}[htbp]
    \centering

    \includegraphics[width=0.48\textwidth]{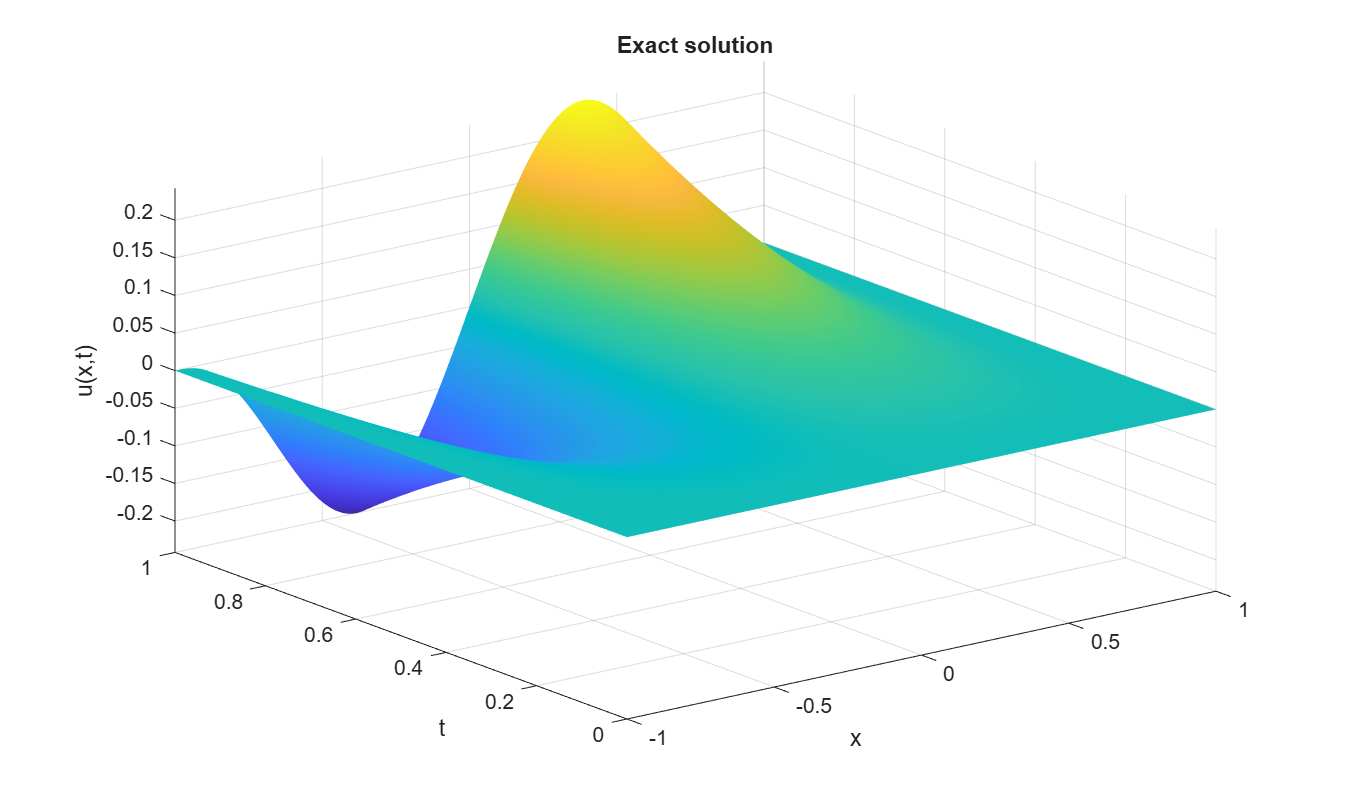}
    \hfill
    \includegraphics[width=0.48\textwidth]{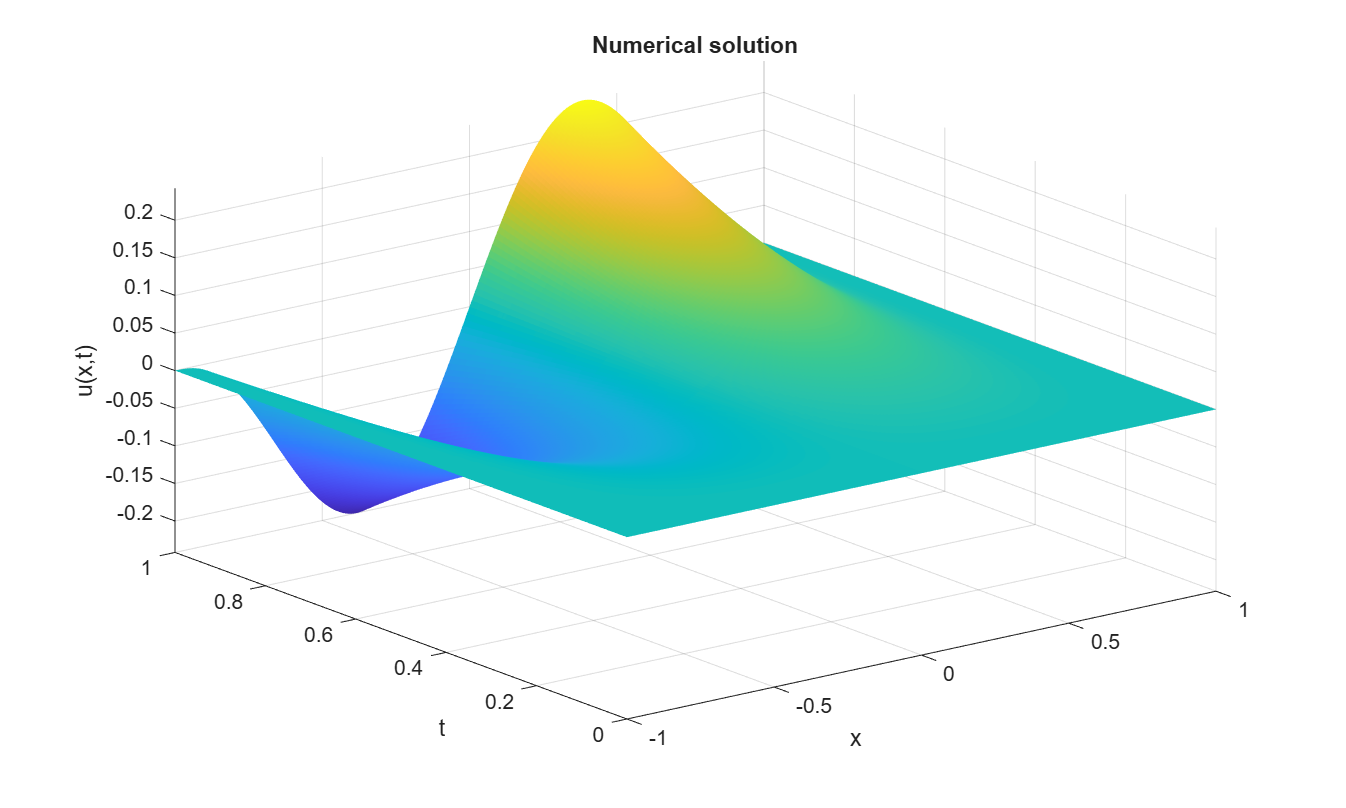}

    \caption{Comparison of exact and numerical solution of Example \ref{Example 6.2} for $\alpha=0.5,\;\beta=1.8,\;N=2048,\;T=1\; \text{and}\; r=(2-\alpha)/\alpha$ .}
    \label{figure2}
\end{figure}

\end{example}
Now, we are considering the domain $(0,1)$ in the later examples to investigate the additional test problems and further examine the applicability of the numerical method to different problem settings.

\begin{example}\label{Example 6.3}
   Consider the test problem with $\Omega\in(0,1),$ 
 \[
 \begin{cases}
       (- \Delta)^{\beta/2} u  = f(x), \quad x \in \Omega,\\
      u(x)  = 0, \quad  x\in \mathbb{R}\setminus\Omega.\\
 \end{cases}
 \]
 With source function \cite{zhang2018riesz}
 \begin{align*}
     f(x)  = 
     \frac{-c_{1,\beta}\Gamma(-\beta)}{\Gamma(4-\beta)} & \left( 2(3-\beta)x^{2-\beta} -6x^{3-\beta}+6(1-x)^{3-\beta} \right. \\
& \left. - 4(3-\beta)(1-x)^{2-\beta}
+(3-\beta)(2-\beta)(1-x)^{1-\beta}
\right),
\end{align*}
for $\beta\in(0,1)\cup (1,2)$. For $\beta=1$, the source term is given as \[
f(x)=\frac {1}{\pi}\left( 3x-\frac{1}{2}+(3x^2-2x) \log \left(\frac{1-x}{x}\right)\right).
\]
 The exact solution $u$ is of the form:
 \begin{align}
 u(x)=
 \begin{cases}
     x^2(1-x),& x\in\Omega,\\
     0,& \text{otherwise.}
 \end{cases}
 \end{align}
  Clearly, the solution is independent of the time variable $t$, so its Caputo fractional derivative with respect to $t$ vanishes. Consequently, the time fractional term does not contribute to the source term $f$.

\begin{table}[htbp]
{ \scriptsize
\centering
\caption{Maximum errors and spatial convergence rate for Example \ref{Example 6.3}.}
\label{Table_4}
\begin{tabular}{c|cc|cc|cc|cc|cc}
\toprule
$(\mu,\beta)$
& \multicolumn{2}{c|}{$(1.8,1.2)$}
& \multicolumn{2}{c|}{$(1.8,1.4)$}
& \multicolumn{2}{c|}{$(1.7,1.5)$}
& \multicolumn{2}{c|}{$(0.8,0.6)$}
& \multicolumn{2}{c}{$(0.8,0.7)$} \\
\cmidrule(lr){2-3}
\cmidrule(lr){4-5}
\cmidrule(lr){6-7}
\cmidrule(lr){8-9}
\cmidrule(lr){10-11}
$M$
& Err & Rate
& Err & Rate
& Err & Rate
& Err & Rate
& Err & Rate \\
\midrule

32
& $4.6328e$-$03$ & --
& $6.6736e$-$03$ & --
& $8.7477e$-$03$ & --
& $4.4525e$-$03$ & --
& $1.2776e$-$02$ & -- \\

64
& $2.7152e$-$03$ & 0.77
& $4.4543e$-$03$ & 0.58
& $6.2711e$-$03$ & 0.48
& $1.7700e$-$03$ & 1.33
& $5.6288e$-$03$ & 1.18 \\

128
& $1.5802e$-$03$ & 0.78
& $2.9564e$-$03$ & 0.59
& $4.4756e$-$03$ & 0.49
& $6.8688e$-$04$ & 1.37
& $2.3795e$-$03$ & 1.24 \\

256
& $9.1551e$-$04$ & 0.79
& $1.9955e$-$03$ & 0.60
& $3.1858e$-$03$ & 0.49
& $2.6341e$-$04$ & 1.38
& $9.8559e$-$04$ & 1.27 \\

512
& $5.2888e$-$04$ & 0.79
& $1.2919e$-$03$ & 0.60
& $2.2635e$-$03$ & 0.49
& $1.0042e$-$04$ & 1.39
& $4.0416e$-$04$ & 1.29 \\

\bottomrule
\end{tabular}
}
\end{table}

\begin{figure}[htbp]
    \centering

    \includegraphics[width=0.60\textwidth]{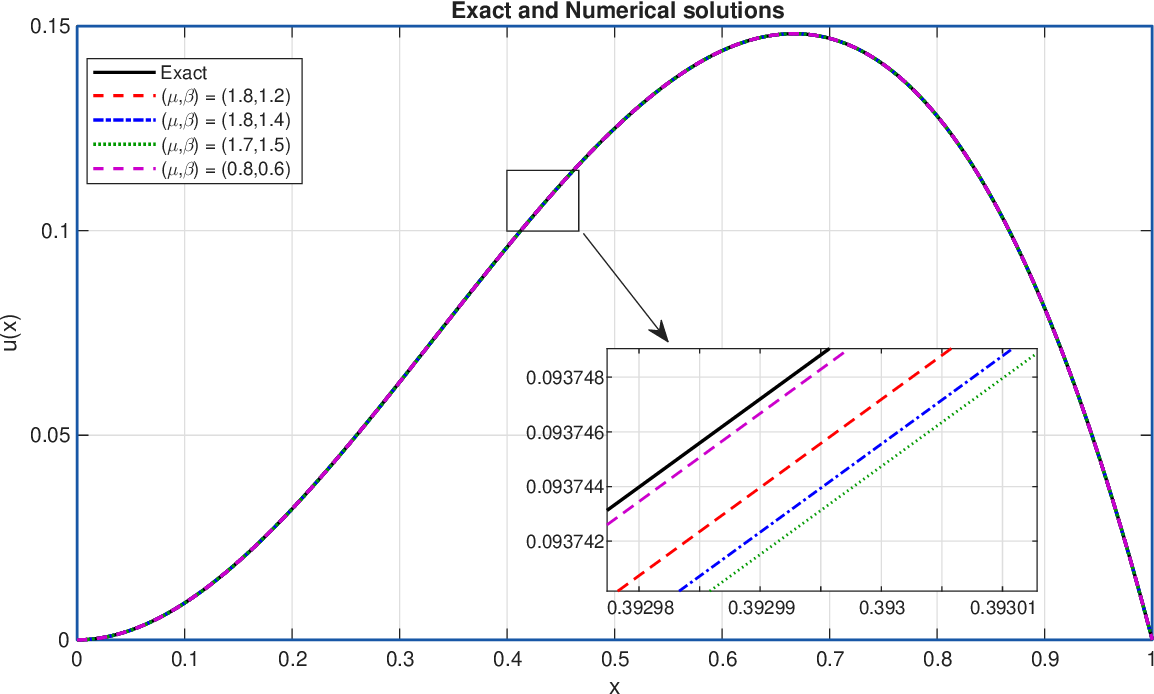}
    \caption{Comparison of the exact and numerical solutions of  Example \ref{Example 6.3} for different values of $\beta$ and $\mu$.}
    \label{figure3}
\end{figure}

From Table \ref{Table_4} and Table \ref{Table_5} we can clearly see that the convergence rate depends upon $\mu$ and $\beta$ continuously. For $\mu\in(\beta,2)$ with $\mu\neq\{2, {1+\beta/2}\}$ and $ \beta\in(1,2)$ the numerical results confirm the theoretical result \eqref{spatial err} of convergence rate $O(h^{2-\beta}).$ Meanwhile, the rate is not observed uniformly  when $\beta\in(0,1]$ for the same $\mu$. However, for $\beta \in (0,1]$, we choose $\mu$ sufficiently close to  $\beta$, which yields convergence rate close to rate $O(h^{2-\beta})$ (see cases of Section \ref{spatial err}). Figure \ref{figure3} and Figure \ref{figure4} we can see that the numerical solution closely agree with the exact solution, causing their overlap. To better visualize this agreement, we have given the enlarged views of selected region of the plots. However, Figure \ref{figure5} exhibit the slower convergence rate, indicating the limited regularity of the exact solution for the proposed numerical problem.

\begin{figure}[htbp]
    \centering
    \includegraphics[width=0.60\textwidth]{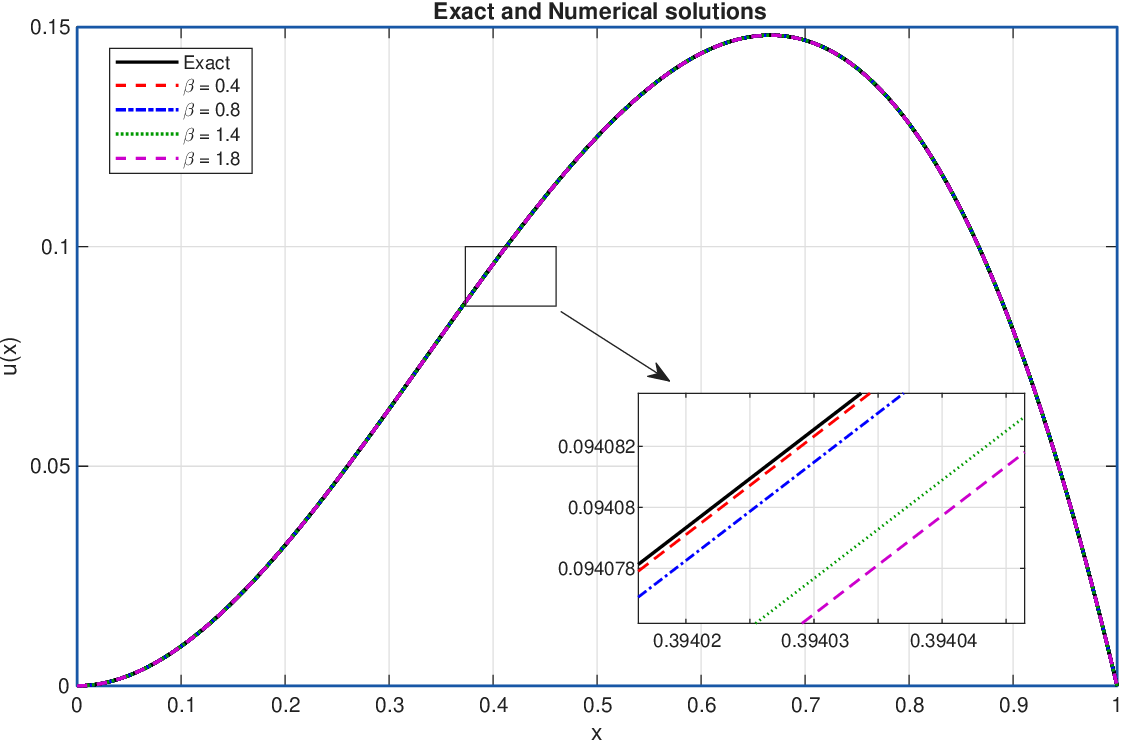}
    \caption{Comparison of the exact and numerical solutions of  Example \ref{Example 6.3} for different values of $\beta$ and $\mu=1+\beta/2$.}
    \label{figure4}
\end{figure}

\begin{figure}[htbp]
    \centering
    \includegraphics[width=0.48\textwidth]{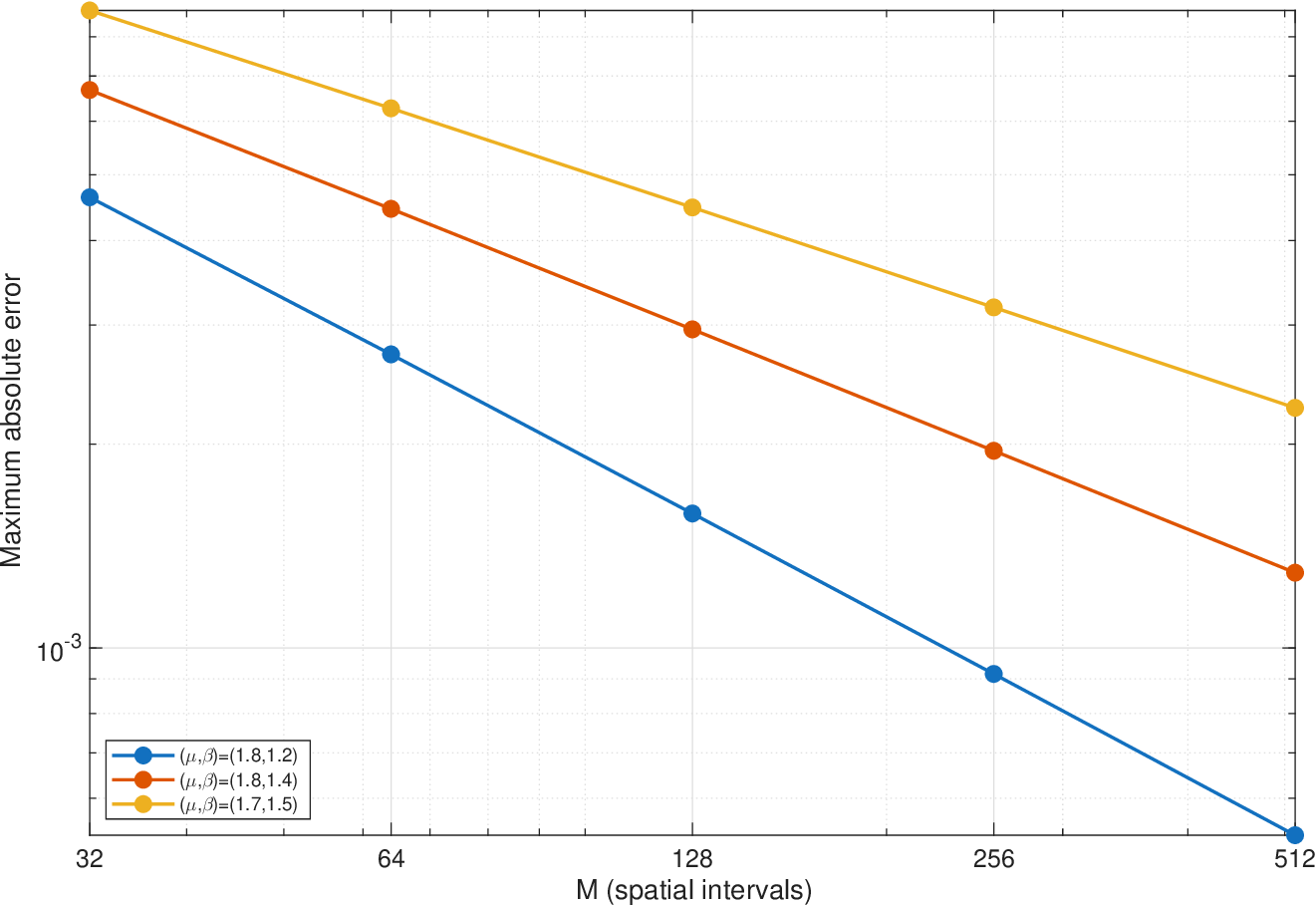}
    \hfill
    \includegraphics[width=0.48\textwidth]{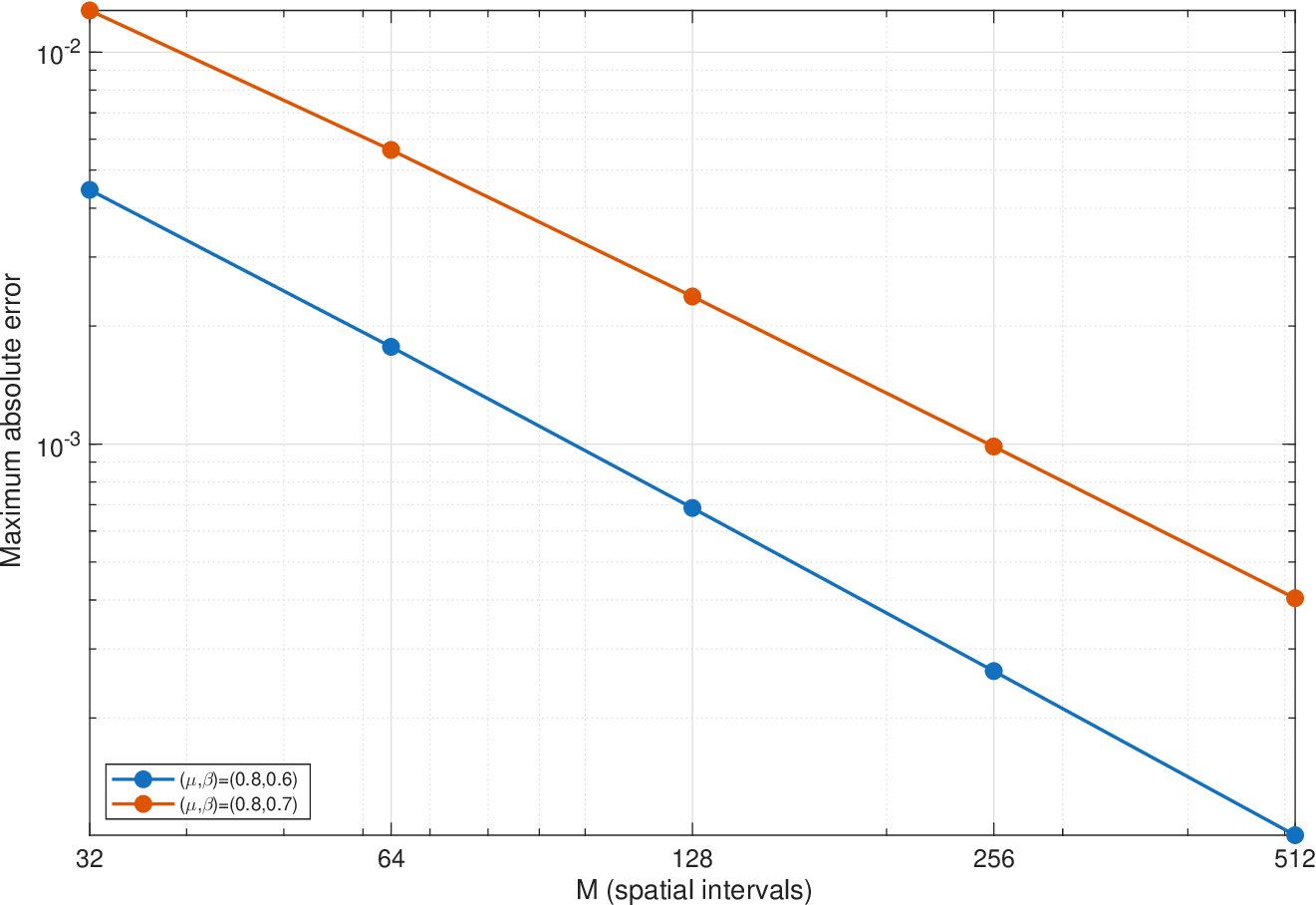}
    \caption{Spatial error decay graph for Table 4 of Example \ref{Example 6.3} .}
    \label{figure5}
\end{figure}
\begin{table}[htbp]
\centering
\caption{Spatial convergence results for different values of  $\beta$ with $\mu=1+\beta/2$ of Example \ref{Example 6.3}.}
\label{Table_5}
\begin{tabular}{c|cc|cc|cc}
\hline
$M$ &
\multicolumn{2}{c|}{$\beta=0.9$} &
\multicolumn{2}{c|}{$\beta=1.4$} &
\multicolumn{2}{c}{$\beta=1.7$} \\
\cline{2-7}
& Error & Rate & Error & Rate & Error & Rate \\
\hline
32 &
$5.6296e$-$04$ & -- &
$5.2720e$-$04$ & -- &
$3.6217e$-$04$ & -- \\

64 &
$2.9492e$-$04$ & $0.93$ &
$2.7519e$-$04$ & $0.94$ &
$1.9276e$-$04$ & $0.91$ \\

128 &
$1.5092e$-$04$ & $0.97$ &
$1.4058e$-$04$ & $0.97$ &
$9.9435e$-$05$ & $0.95$ \\

256 &
$7.6335e$-$05$ & $0.98$ &
$7.1048e$-$05$ & $0.98$ &
$5.0617e$-$05$ & $0.97$ \\

512 &
$3.8388e$-$05$ & $0.99$ &
$3.5715e$-$05$ & $0.99$ &
$2.5565e$-$05$ & $0.99$ \\

\bottomrule
\end{tabular}

\end{table}
It is observed that the exact solution $u(x)=x^2(1-x)$  satisfies the prescribed boundary condition at both the end points. However $u'(x)=2x-3x^2,$  gives $u'(0)=0 \;\text{and}\; u'(1)=-1\neq0$, hence $u$ vanishes to the second order at $x=0,$ but first order at $x=1$, shows that the first order derivative has a jump discontinuity at $x=1$. Consequently, the zero extension of $u$ outside the computational domain $(0,1),$ which result the function continuous across $x=1$, doesn't possess the higher order regularity required for the optimal second-order $O(h^2)$ spatial approximation (shown in Table \ref{Table_5}). Therefore despite of the fact that the exact solution satisfies the prescribed boundary condition, lack of sufficient regularity at $x=1,$ limits the spatial convergence rate to $O(h).$ In particular, theoretical second-order convergence rate is not attained for $\mu=2.$

 \end{example}

 Now, we consider an example where the source term is absent and the exact solution is unknown. To investigate the convergence behavior of the numerical solution, we employ the two-mesh principle, discussed at the beginning of Section 6.

\begin{example}\label{Example 6.4}
 We consider the following test problem:
 \[
 \begin{cases}
      _{}^CD_t^\alpha u + (- \Delta)^{\beta/2} u  = 0, \quad (x,t) \in (0,1) \times (0,1],\\
      u(x,t)  = 0, \quad x\in\mathbb{R}\setminus(0,1),\;t \in (0,1],\\
      u(x,0)  = \sin(\pi x), \quad  x \in (0,1).
 \end{cases}
 \]
 We apply the numerical scheme with choosing $M$ and $N$ with different choices of $\alpha$ and $\beta$.

\begin{figure}[htbp]
    \centering
    \includegraphics[width=0.7\textwidth]{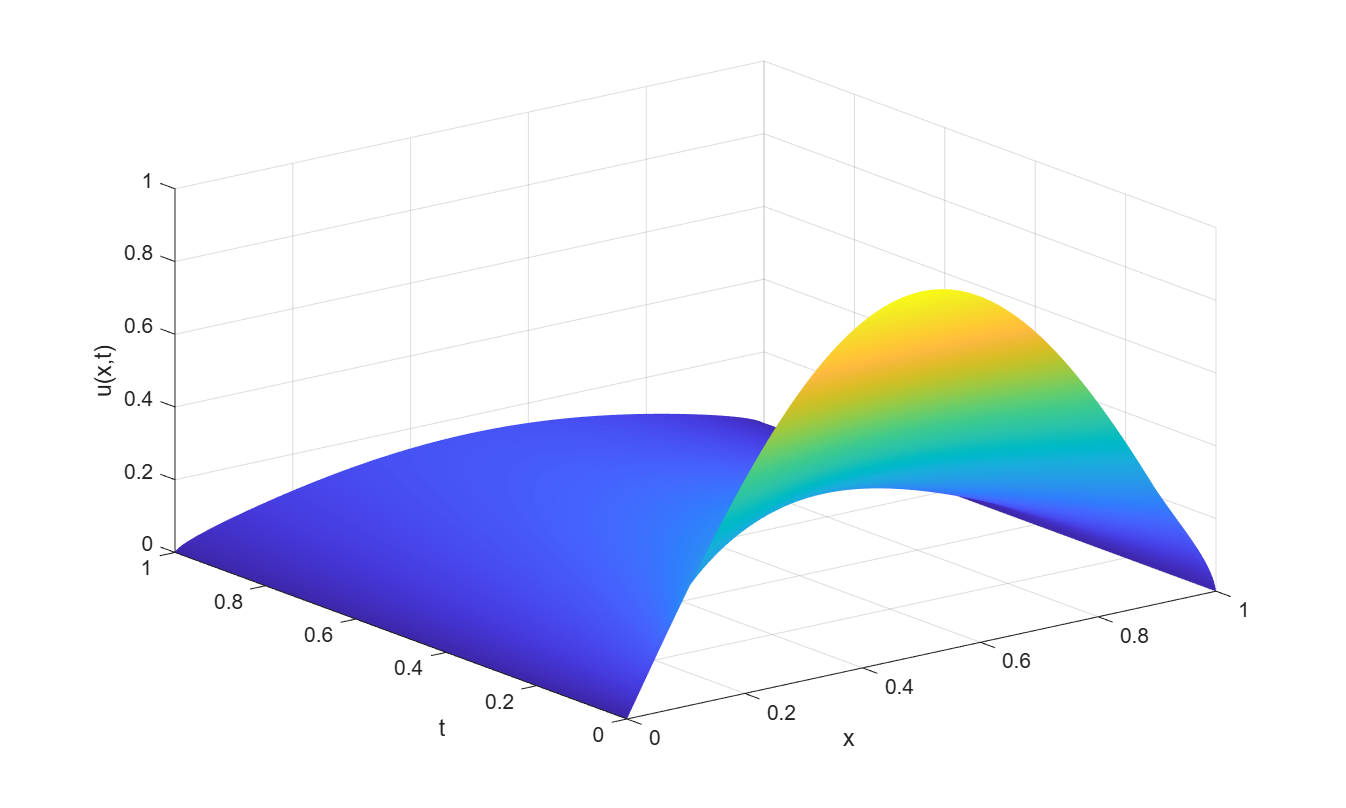}
    \caption{Numerical solution for Example \ref{Example 6.4} when $\alpha=0.6, \beta=1.2$ at $N=1024$, $M=32$.}
    \label{figure6}
\end{figure}

 \begin{table}[htbp]
\centering
\caption{Spatial convergence results for different values of $\alpha$ and $\beta$ with $\mu$=$1+\beta/2$ and $N=1024$ for Example \ref{Example 6.4}.}
\label{Table_6}
\begin{tabular}{c|cc|cc|cc|cc}
\hline
$M$ &
\multicolumn{2}{c|}{$(\alpha,\beta)=(0.8,1.7)$} &
\multicolumn{2}{c|}{$(\alpha,\beta)=(0.4,1)$} &
\multicolumn{2}{c|}{$(\alpha,\beta)=(0.6,1.2)$} &
\multicolumn{2}{c}{$(\alpha,\beta)=(0.7,1.5)$}\\
\cline{2-9}
& Error & Rate & Error & Rate & Error & Rate  & Error & Rate \\
\hline
32 &
 -- & -- &
-- & -- &
-- & -- &
-- & -- \\

64 &
$2.45061e$-$03$ & -- &
$5.86447e$-$03$ & -- &
$5.65224e$-$03$ & -- &
$1.25769e$-$03$ & --\\

128 &
$1.51655e$-$03$ & $0.69$ &
$4.58490e$-$03$ & $0.36$ &
$4.11239e$-$03$ & $0.46$ &
$7.56550e$-$04$ & $0.73$\\
256 &
$8.98740e$-$04$ & $0.75$ &
$3.47432e$-$03$ & $0.40$ &
$2.88876e$-$03$ & $0.51$ &
$4.52105e$-$04$ & $0.74$\\
512 &
$5.19974e$-$04$ & $0.79$ &
$2.57226e$-$03$ & $0.43$ &
$1.98101e$-$03$ & $0.54$ &
$2.69460e$-$04$ & $0.75$\\
1024 &
$2.99229e$-$04$ & $0.80$ &
$1.87308e$-$03$ & $0.47$ &
$1.33722e$-$03$ & $0.57$ &
$1.60408e$-$04$ & $0.75$\\
\bottomrule
\end{tabular}

\end{table}
From the Table \ref{Table_6}, we observe that the spatial convergence rate is approximately $O(h^{\beta/2})$ and Table \ref{Table_7} indicates a first-order temporal convergence rate. However, both the observed rates are suboptimal compared with the expected second-order spatial rate and theoretical temporal rate of $O(N^{-(2-\alpha)}).$
Since the exact solution is unknown, the reduced convergence rates may be attributed to limited regularity of the solution, which can affect both the space and time discretizations.
\begin{table}[htbp]
\centering
\caption{Temporal convergence results for different values of $\alpha$ and $\beta$ with $\mu$=1+$\beta/2$, $r=(2-\alpha)/\alpha$ \; and $M=1024$ for Example \ref{Example 6.4}}
\label{Table_7}
\begin{tabular}{c|cc|cc|cc}
\hline
$N$ &
\multicolumn{2}{c|}{$(\alpha,\beta)=(0.3,0.6)$} &
\multicolumn{2}{c|}{$(\alpha,\beta)=(0.2,0.4)$} &
\multicolumn{2}{c}{$(\alpha,\beta)=(0.4,0.8)$} \\
\cline{2-7}
& Error & Rate & Error & Rate & Error & Rate \\
\hline
32 &
-- & -- &
-- & -- &
-- & -- \\
64 &
$2.89528e$-$02$ & -- &
$5.10975e$-$02$ & -- &
$1.46517e$-$02$ & -- \\

128 &
$2.26067e$-$02$ & $0.36$ &
$4.02689e$-$02$ & $0.34$ &
$1.14228e$-$02$ & $0.36$ \\

256 &
$1.45375e$-$02$ & $0.64$ &
$2.55587e$-$02$ & $0.66$ &
$7.25449e$-$03$ & $0.65$ \\

512 &
$8.22346e$-$03$ & $0.82$ &
$1.43361e$-$02$ & $0.83$ &
$3.86391e$-$03$ & $0.91$ \\

1024 &
$4.27202e$-$03$ & $0.94$ &
$7.50512e$-$03$ & $0.93$ &
$1.75027e$-$03$ & $1.14$ \\

2048 &
$2.08458e$-$03$ & $1.03$ &
$3.77727e$-$03$ & $0.99$ &
$7.33952e$-$04$ & $1.24$ \\

\bottomrule
\end{tabular}

\end{table}
\end{example}

\section{\textbf{Conclusion }}\label{Sec.Conclude}
In this paper, we studied a Dirichlet boundary value problem for a time-fractional diffusion equation involving fractional Laplacian operator. We established a suitable weak formulation for the problem within the framework of fractional Sobolev space. Under the appropriate assumption on source term and initial data, the existence and uniqueness of the weak solution are proved using \vspace{0.5cm}standard tools from functional analysis and some prior results.\\
For the numerical approximation Caputo derivative was discretized using temporal graded mesh using $L1$ finite difference scheme, while the spatial discretization is done by using appropriate fractional Laplacian approximation.  Convergence of the proposed scheme was established by proving the error analysis.
The theoretical result demonstrate that the method is reliable for \vspace{0.5cm}solving a class of fully non-local diffusion problem.\\
Finally, numerical experiments are conducted to validate the theoretical error estimates, confirming the accuracy and robustness of the scheme.\\
Future work may focus on extending the present analysis on higher dimensional domain, developing more efficient numerical algorithms for the approximation of fractional Laplacian, different boundary conditions, etc.


\end{document}